\documentclass[twoside,leqno]{article}

\usepackage[letterpaper]{geometry}

\usepackage{siamproceedings}
\usepackage[utf8]{inputenc}
\usepackage{color}
\usepackage{makeidx}
\usepackage{xspace}
\usepackage{amsmath}
\usepackage{amssymb}
\usepackage{amsfonts}
\newsiamremark{remark}{Remark}

\newcommand{\dP}{\mathbb{P}}

\newcommand{\dN}{\mathbb {N}}
\newcommand{\dR}{\mathbb {R}}

\newcommand{\cF}{\mathcal {F}}

\newcommand{\cB}{\mathcal {B}}
\newcommand{\cC}{\mathcal {C}}
\newcommand{\cG}{\mathcal {G}}

\newcommand{\cX}{\mathcal {X}}

\newcommand{\cP}{\mathcal {P}}
\newcommand{\cp}{\mathfrak {p}}
\newcommand{\cQ}{\mathcal {Q}}

\newcommand{\cI}{\mathcal {I}}
\newcommand{\cN}{\mathcal {N}}

\newcommand{\cS}{\mathcal {S}}

\newcommand{\cT}{\mathcal {T}}
\newcommand{\ct}{\mathfrak{t}}
\newcommand{\cY}{\mathcal {Y}}

\newcommand{\II}{1\!\!{\sf I}}

\newlength{\figurewidth}
\title{Computational thresholds in high-dimensional statistics: the case of graph alignment}
\author{Laurent Massoulié\thanks{Inria, DI-ENS PSL Research University}}
\begin{document}
\maketitle
\begin{abstract} 

Graph alignment is an instance of the NP-hard quadratic assignment problem. It consists, given the adjacency matrices $A_1$, $A_2$ in $\dR^{n\times n}$ of two graphs, in finding a permutation matrix $\Pi$  that minimizes the Frobenius norm $\| \Pi A_1-A_2 \Pi\|_F$. In this article we consider the problem from the perspective of high-dimensional statistics: we aim to estimate an unknown permutation $\pi^*$ in the symmetric group $\cS_n$ from the observation of two correlated random adjacency matrices $A_1$, $A_2$.

\hbox{}

We establish the following computational thresholds. We first consider the case where $A_1$, $A_2$ are the adjacency matrices of two correlated Erd\H{o}s-Rényi random graphs $\cG(n,p)$ in the sparse regime with average degree $\lambda:=np= O(1)$ and edge correlation parameter $s\in(0,1)$. We identify a critical threshold $s^*(\lambda)$ for $s$ above which a message-passing, local algorithm succeeds at alignment --that is, recovers a fraction $\Omega(1)$ of the entries of $\pi^*$, and below which no local algorithm succeeds. This result crucially depends on an associated model of correlated random trees. 
%We review existing conjectural evidence that when this algorithm fails, no other polynomial-time algorithm can succeed.

\hbox{}

We then consider the case where $A_1$, $A_2$ are two correlated Gaussian Wigner matrices with correlation parameter expressed as $s=1/\sqrt{1+\sigma^2}$ where $\sigma$ is the ``noise'' parameter. In this setting we consider a fast spectral algorithm based on the leading eigenvectors of the matrices $A_1$, $A_2$, and identify the critical scaling for noise parameter $\sigma$ at which the fraction of entries of $\pi^*$ correctly recovered goes from $1-o(1)$ to $o(1)$. We next consider the convex relaxation approach which obtains the doubly stochastic matrix $X$ that minimizes $\|X A_1 -A_2 X\|_F$. We obtain upper and lower bounds on the critical noise parameter $\sigma$ at which a simple post-processing of $X$ correctly recovers a fraction $1-o(1)$ of entries of $\pi^*$.

\hbox{}

We finally identify promising future directions on i) computational thresholds for spectral methods and convex relaxation methods of practical interest, and ii) impossibility results for broad classes of algorithms, notably low degree polynomial algorithms and local search algorithms.

%We next discuss a popular convex relaxation which returns the doubly stochastic matrix $X$ that minimizes $\|X A_1 - A_2 X\|_F$. We show that, for correlation parameter $1/\sqrt{1+\sigma^2}$ between the two matrices, and any $\epsilon>0$, denoting by $\Pi^*$ the target permutation matrix, with high probability $\|\Pi^*-X\|^2_F=o(n)$ if $\sigma=n^{-1/2-\epsilon }$ and $\|\Pi^*-X\|^2_F=\Omega(n)$ if instead $\sigma=n^{-1/2+\epsilon }$. Thus at threshold $\sigma=1/\sqrt{n}$, algorithms based on this convex relaxation undergo a qualitative change in performance.  Promising future directions are finally given.
\end{abstract}

\section{Introduction.} 
High-dimensional statistics focuses on statistical estimation problems where the observed data, that we shall denote by $A$, and the signal to be estimated, that we shall denote by $\pi^*$, are both ``high-dimensional''. A prevalent phenomenon in this context is the so-called ``informational-computational gap'': For many problems of interest, there exist parameter ranges where brute-force search yields some non-trivial estimator $\hat \pi$, while all known polynomial-time algorithms fail to do so. 
%This can be understood in terms of phase diagrams: for high noise levels, estimation is impossible for lack of sufficient information in the observations; for low noise levels, estimation can be done successfully with adequate polynomial-time algorithms; for intermediate noise levels, brute force succeeds while known polynomial-time methods fail.

\hbox{}

A major objective of the field is to identify efficient polynomial-time algorithms and determine under which conditions they give non-trivial estimators $\hat \pi$. Convex relaxations, spectral methods, message-passing algorithms and local search methods such as Monte-Carlo Markov Chains (MCMC) are the natural candidates in this endeavour.

\hbox{}

An equally important objective is to determine conditions under which a whole family of algorithms must fail. One approach consists in studying  the ``energy landscape'' $\pi\to \dP(\pi^*=\pi|A)$, assuming a natural topology on the space $\cS$ to which $\pi^*$ belongs. This may allow one to rule out local search algorithms that search for local maxima of the posterior distribution $\pi\to \dP(\pi^*=\pi|A)$. One then says that the corresponding problem is ``local search-hard''. Bandeira et al. \cite{franz-parisi} provide criteria for this failure mode of local search algorithms based on constructs from statistical physics, namely the so-called Franz-Parisi potential.

\hbox{}

Another fruitful approach consists of ruling out estimates $\hat\pi$ that depend polynomially on the observation $A$, the degree of the polynomial being bounded by $D=O(\ln n)$ where $n$ captures the problem's dimension; see Wein \cite{Wein-arxiv} for a recent survey and Hopkins \cite{thesis-hopkins} for an early reference. One then says that the problem at hand is ``LDP-hard'', where LDP stands for Low Degree Polynomial. 

\hbox{}

Correspondence between the two notions of hardness was investigated in Bandeira et al. \cite{franz-parisi}, where they showed that the two types of hardness coincide for problems of detecting whether some signal has been added or not to a vector of Gaussian noise. 

\hbox{}

Together, positive results on some algorithm and negative results on a whole family may combine to show that this specific algorithm is optimal within that family. 

\hbox{}

Our goal is to explain results on the efficiency of illustrative algorithms and the underlying mathematical techniques in the context of graph alignment. The graph alignment problem is popular because it arises in a variety of applications, e.g. network de-anonymization (see Narayanan and Shmatikov \cite{narayanan2008robust}), computational biology (see Singh et al. \cite{singh2008global}), pattern recognition (see Conte et al. \cite{conte2004thirty}) to name a few. It has also attracted attention in high-dimensional statistics recently as its analysis required the development of new techniques.
%requires new techniques beyond those developed for the study of better-studied problems such as graph clustering.  

\hbox{}

In graph alignment, one observes two adjacency matrices $A_1,\; A_2\in \dR^{n\times n}$, generated in a two-step procedure. In a first step, $A_1$ is generated at random and a noisy version $A'_2=A_1+W$ of $A_1$ is obtained through addition of noise matrix $W$. In a second step, a permutation $\pi^*$ is sampled uniformly at random from the symmetric group $\cS_n$, independently of $(A_1,A'_2)$, and matrix $A_2$ is then defined via $A_2(i,j)=A'_2(\pi^*(i),\pi^*(j))$.

\hbox{}

Consistent with previous notation, $\pi^*$ is the signal of interest, to be estimated
%Our goal is then to estimate the signal of interest, that is the permutation $\sigma^*$, 
from observation $A=(A_1,A_2)$. The quality of estimator $\hat\pi\in \cS_n$ is assessed by its so-called overlap 
$$
\hbox{ov}(\pi^*,\hat\pi):=\frac{1}{n}\sum_{i\in[n]}\II_{\pi^*(i)=\hat\pi(i)},
$$ 
that is the fraction of correctly estimated entries of $\pi^*$.
Several objectives could be envisioned, expressed in terms of $\hbox{ov}(\pi^*,\hat\pi)$: exact reconstruction corresponds to $\hbox{ov}(\pi^*,\hat\pi)=1$ and almost exact reconstruction to $\hbox{ov}(\pi^*,\hat\pi)=1-o(1)$.
Of particular interest is partial reconstruction, that is when when $\hbox{ov}(\pi^*,\hat\pi)\ge\Omega(1)$ as $n\to\infty$. 

\hbox{}

Graph alignment as described is an instance of Bayesian estimation: $\pi^*$ is the signal --or structure-- planted in the observation, sampled from a prior distribution, that we aim to estimate. 
As is typical in high-dimensional statistics, one can associate to this estimation problem a related hypothesis testing, or decision problem: the couple $(A_1,A_2)$ generated as just described depends on signal $\pi^*$; alternatively one could have generated $A_1,A_2$ independently of one another, but with the same marginal distributions. The latter situation would be the null hypothesis $H_0$ (no signal), and the former the alternative hypothesis $H_1$ (signal present). The decision problem then consists in testing from sample $(A_1,A_2)$ whether one is under $H_0$ or $H_1$. This correspondence between a reconstruction and a decision problem is reminiscent of computational complexity theory where optimization problems are also associated with decision problems.

\hbox{}

%To make the discussion more concrete, let us precise two specific versions of the problem.

The two most studied versions of the problem are the following:

\hbox{}

{\bf Gaussian ensembles}: $A_1$ is a Gaussian Wigner matrix, also known as the Gaussian Orthogonal Ensemble (GOE), that is: $(A_1(i,j))_{i\le j\in [n]}$ are mutually independent, $A_1(i,j)=A_1(j,i)$, $A_1(i,i)\sim \cN(0,2/\sqrt{n})$, $i\in[n]$, and $A_1(i,j)\sim \cN(0,1/\sqrt{n})$, $i<j\in [n]$. Matrix $A'_2$ reads $A'_2=A_1+ \sigma Z$, where $Z$ is independent of $A_1$ and identically distributed to $A_1$ and $\sigma$ is the noise parameter.

\hbox{}

{\bf Erd\H{o}s-Rényi random graphs}: here $A_1$ is the adjacency matrix of an Erd\H{o}s-Rényi random graph, i.e. $(A_1(i,j))_{i<j}$ are i.i.d. Bernoulli$(p)$ for some $p\in [0,1]$, $A_1(i,i)\equiv 0$, and $A_1(i,j)=A_1(j,i)$, $i<j\in [n]$. For a given entry $(i,j)$, $i<j$, the pair $(A_1(i,j),A'_2(i,j))$ follows the Bernoulli$(p,s)$ distribution defined via
$$
(A_1(i,j),A'_2(i,j))=\left\{\begin{array}{ll}(1,1)&\hbox{with probability }p s,\\
(1,0)& \hbox{with probability }p (1-s),\\
(0,1)& \hbox{with probability }p (1-s),\\
(0,0)& \hbox{with probability }1-p(2-s).
 \end{array}
\right.
$$
These two versions typically rely on different sets of tools; random matrix theory for the former, and random graph theory for the latter. 

\hbox{}

\subsection{Outline:}
In {\bf Section \ref{sec:2}} we focus on the Erd\H{o}s-Rényi model in the sparse regime where the average degree $\lambda:=np$ is of constant order $O(1)$. This regime is of special interest for several reasons: i) it is challenging for the design of efficient algorithms; ii) sparse graphs occur in application scenarios; iii) its analysis can be done by the study of well-defined limiting local processes, random trees as it will appear.

\hbox{}

We introduce MP-Align, a message-passing algorithm which matches vertex $i$ in $A_1$ to vertex $u$ in $A_2$ if their local neighborhoods in their respective graphs are ``sufficiently similar''. We leverage the local structure of Erd\H{o}s-Rényi graphs in the sparse regime to characterize the critical threshold $s^*(\lambda)$ for $s$ above which MP-Align succeeds. We also establish that for $s$ below this threshold, for a suitable notion of local algorithm, no $O(\ln (n))$-local algorithm succeeds. 

\hbox{}

In {\bf Section \ref{sec:3}} we relate the success of MP-ALign to properties of a hypothesis testing problem, deciding whether a pair of random trees is correlated or independent. We provide an asymptotic characterization of this tree correlation testing problem in the limit $\lambda\to\infty$, establishing in particular that $\lim_{\lambda\to\infty}s^*(\lambda)=\sqrt{\alpha}$, where $\alpha$ is Otter's constant \cite{otter1948number}. 
%We conclude the section by discussing related results in the literature, motivating the conjecture that below the threshold $s^*(\lambda)$, no polynomial-time algorithm succeeds.

\hbox{}

{\bf Section \ref{sec:spectral}} is devoted to the study of a specific fast spectral algorithm that relies on the leading eigenvectors of both matrices $A_1$, $A_2$. There we place ourselves in the Gaussian ensemble setting. We establish that for any $\epsilon>0$, the algorithm achieves overlap $1-o(1)$ for noise parameter $\sigma=n^{-7/6-\epsilon}$, and instead achieves overlap $o(1)$ for noise parameter $\sigma=n^{-7/6+\epsilon}$. 

\hbox{}

{\bf Section \ref{sec:4}} is devoted to the convex relaxation where rather than minimizing $\|X A_1 -A_2 X\|_F$ over permutation matrices $X$, one relaxes this to minimization over the convex hull of permutation matrices, also known as the Birkhoff polytope, or equivalently the set of doubly stochastic matrices. Despite this relaxation being very successful and popular in practice, there are few rigorous results to date on their performance.
%as their analysis is very challenging. 
 Focusing on the Gaussian ensemble case, we prove that for any $\epsilon>0$, denoting by $\Pi^*$ the optimal permutation matrix, with high probability $\|\Pi^*-X\|^2_F=o(n)$ if $\sigma=n^{-1-\epsilon }$ and $\|\Pi^*-X\|^2_F=\Omega(n)$ if instead $\sigma=n^{-1/2+\epsilon }$. Thus algorithms based on this convex relaxation undergo a qualitative change in performance for $\sigma$ between $n^{-1}$ and $n^{-1/2}$. 

\hbox{}

{\bf Section \ref{sec:5}} concludes by reviewing promising future directions. Specifically, Sections \ref{sec:2} and \ref{sec:3} lead to conjectures on LDP-hardness as well as local search-hardness of graph alignment for sparse correlated Erd\H{o}s-Rényi graphs. Sections \ref{sec:spectral} and \ref{sec:4} lead to several questions respectively on the spectra of composite random matrices, and on finer performance guarantees for convex relaxations.
%alternative spectral methods. 

\hbox{}

\subsection{Related works:}
Study of graph alignment for Erd\H{o}s-Rényi random graphs initially focused on informational thresholds, i.e. feasibility of the problem irrespective of computational complexity. Cullina and Kyavash \cite{cullina-kyavash} showed that exact recovery of unknown permutation $\pi^*$ is feasible with high probability if and only  $ps-\ln(n)=\omega(1)$, or equivalently if and only if the intersection graph (whose adjacency matrix $A_\cap$ is given by $A_1\odot A'_2$) is fully connected with high probability. The study of partial recovery (achieving overlap $\Omega(1)$) started with Hall and Massoulié \cite{DBLP:journals/ior/HallM23}; partial recovery is feasible if and only if $s np>1$ or equivalently if and only if the intersection graph admits a giant (i.e. comprising $\Omega(n)$ vertices) connected component. The necessary part was established by  Ganassali, Lelarge and Massoulié \cite{DBLP:conf/colt/GanassaliML21} and the sufficient part by Ding and Du \cite{ding-hu}. 

\hbox{}

A series of works proposed polynomial-time algorithms and proofs of their success, for either exact or partial recovery of $\pi^*$. In particular Ding et al. 
\cite{ding2021efficient} obtain exact recovery for $np\ge \ln^{\alpha}(n)$ and $1-s\le \ln^{-\beta}(n)$ for suitable constants $\alpha,\; \beta>0$; Mao et al. \cite{mao2023exact} relaxed the condition on $s$ to $1-s\le c$ for some constant $c\in (0,1)$. It was then shown by Mao et al. \cite{mao2023random} that partial recovery can be done in polynomial-time for graphs with average degree $np$ sufficiently large, as soon as $s>\sqrt{\alpha}$, where $\alpha$ is Otter's constant.

\hbox{}

Concerning impossibility results, Ding, Lu and Li \cite{ding-du-li-23} showed that the hypothesis testing problem of distinguishing between correlated graphs and independent graphs cannot be done with low degree polynomials whenever $s<\sqrt{\alpha}$. Li \cite{zhangsong-li-25} then showed that partial recovery of $\pi^*$ is impossible with low degree polynomials when $s<\sqrt{\alpha}$.

\hbox{}

For correlated Gaussian Wigner matrices, 
 Ganassali \cite{DBLP:conf/msml/Ganassali21} established that the threshold for informational feasibility of recovery is for correlation parameter $1/\sqrt{1+\sigma^2}$ between the two matrices equal to $2\sqrt{\ln n/n}$. A spectral polynomial-time algorithm, coined GRAMPA, was proposed by Fan et al. \cite{DBLP:journals/focm/FanMWX23}, achieving recovery for $\sigma=O(1/\ln(n))$. The authors extended the results for GRAMPA to the setting of Erd\H{o}s-Rényi graphs in \cite{fan2023spectral2}. Ding and Li \cite{DBLP:journals/focm/DingL25} then obtained a polynomial-time algorithm that succeeds at recovery for any non-vanishing correlation parameter. 

\hbox{}

Much less is known about the performance of convex relaxation approaches in the settings we consider. Araya and Tyagi \cite{araya2024graph} consider convex relaxation from permutation matrices to both the Birkhoff polytope $\cB_n$ and the simplex $\Delta_n:=\{M\in \dR_+^{n^2}: \sum_{i,j\in [n]}M_{ij}=n\}$ in the noiseless case $\sigma=0$, and show that these relaxations then recover $\pi^*$ exactly. 

\hbox{}

%In \cite{sushil25} we show that relaxation to the Birkhoff polytope  yields a doubly stochastic matrix $X^*$ such that its error from the sought permutation matrix $\Pi^*$ verifies $\|\Pi^*-X^*\|_F^2=o(n)$ for $\sigma=n^{-1-\epsilon}$, and  $\|\Pi^*-X^*\|_F^2=\Omega(n)$ for $\sigma=n^{-1/2+\epsilon}$. Our goal in Section \ref{sec:4} is precisely to tighten the gap between these two conditions. 

The theorems stated in Sections \ref{sec:2}  and \ref{sec:3} stem from Ganassali, Lelarge and Massoulié \cite{luca} and Ganassali, Massoulié and Semerjian \cite{ganassali2022statistical} respectively, and we refer to these for proofs (see also Maier and Massoulié \cite{jakob} for detailed proofs, and generalization of the results to a model of asymmetric random graphs with differing numbers of vertices and densities of edges). The results of Section \ref{sec:spectral} appeared in Ganassali, Lelarge and Massoulié \cite{ganassali2022spectral}, to which we refer for detailed proofs. The result of Section \ref{sec:4} is from Varma, Waldspurger and Massoulié \cite{sushil25} and we refer to that paper for its proof.

\hbox{}

% improve known performance bounds for its resolution based on low-degree polynomials. Show that in sparse regime, the precise threshold below which ``local algorithms'' fail is also the threshold below which low-degree polynomials fail.

\section{Sparse Erd\H{o}s-Rényi graphs and a message-passing algorithm.}  \label{sec:2}

To motivate the algorithm, consider the so-called union graph $A_\cup$ of $A_1$ and $A'_2$, where we let 
$$
A_\cup(i,j)=\left\{\begin{array}{ll}\{1,2\}&\hbox{ if }A_1(i,j)=A'_2(i,j)=1;\\
\{1\}&\hbox{if }A_1(i,j)=1,\;A'_2(i,j)=0;\\
\{2\}&\hbox{if }A_1(i,j)=0,\;A'_2(i,j)=1;\\
\emptyset&\hbox{if }A_1(i,j)=A'_2(i,j)=0.
\end{array}\right.
$$
A classical result of random graph theory  states that for an Erd\H{o}s-Rényi random graph $\cG(n,p)$ with $\lambda=np$ of order 1, the neighborhood of a fixed vertex $i$ of the graph up to distance $d=O(\ln n)$ converges, for the total variation distance, to a Galton-Watson branching tree of depth $d$ and with offspring distribution following the Poisson distribution with parameter $\lambda$, that we denote by $\cp_\lambda$. See for instance Alon and Spencer \cite{alon2016probabilistic}. Writing $d=\kappa \ln(n)$, the constant  $\kappa$ ought to be chosen such that $\lambda^d\ll \sqrt{n}$, i.e. $\kappa\ln(\lambda)<1/2$ for this property to hold. 

\hbox{}

This result easily generalizes to a situation where edges have multiple types such as the model of the union graph $A_\cup$ we just defined, where edge types are $\{1\}$, $\{2\}$ or $\{1,2\}$, an entry of $A_\cup$ equal to $\emptyset$ corresponding to absence of an edge; see Bollob{\'a}s Janson and Riordan \cite{bollobas2007phase}. Up to distance $d\le \kappa \ln n$, where $\kappa \ln(\lambda(2-s))<1/2$, the neighborhood of a vertex in the union graph converges in total variation distance to a multi-type Galton-Watson branching process with edge types $\{1\}$, $\{2\}$ or $\{1,2\}$. In this multivariate version, the distribution of the number $(k_1,k_2,k_{12})$ of children of respective types $\{1\}$, $\{2\}$, $\{1,2\}$ (where the type of a child corresponds to the edge between itself and its parent), for any given individual follows the product distribution $\cp_{\lambda(1-s)} \otimes \cp_{\lambda(1-s)}\otimes\cp_{\lambda s}$. 

\hbox{}

We shall use the following           
\begin{definition}
Denote by $\cY_d$ the collection of rooted, labeled trees of depth at most $d$, where the root label is $\bullet$ and for any vertex with label $i$ and $k$ children, its children are labeled as $(i;1),\ldots,(i;k)$.

Denote by $GW_d^{(\lambda)}$ the distribution over $\cY_d$ of a depth $d$-Galton-Watson tree with offspring distribution $\cp_\lambda$.
\end{definition}

\hbox{}

Let us now define a distribution over correlated pairs of Galton-Watson trees as follows:
\begin{definition}
Let parameters $\lambda>0$, $s\in (0,1)$ and depth $d\in \dN$ be fixed. The $(\lambda,s,d)$ {\bf augmentation} of any tree $t\in \cY_d$ is defined as follows: first, for each $h=0,\ldots,d-1$,  endow each depth $h$-vertex of $t$ with an additional number of children distributed according to $\cp_{\lambda(1-s)}$. Next, to each such new vertex, whose depth is $h+1$, graft a descendance tree of depth $d-h-1$ sampled according to $GW^{(\lambda)}_{d-h-1}$, all child additions and tree graftings being performed independently.

\hbox{}

The {\bf shuffling} of any tree $t$ in $\cY_d$ is obtained by successively performing independent uniform permutations of the labels of the children of one vertex, starting from the root and proceeding towards the leaves. The result is a tree from $\cY_d$ sampled uniformly at random from those trees $t'$ that are isomorphic to $t$, for the natural notion of graph isomorphism between rooted labeled trees. 

\hbox{}

$\cP^{(\lambda,s)}_d$ is then defined as the distribution over $\cY_d\times \cY_d$ of a pair $(\cT_d,\cT'_d)$ obtained as follows.  Start from some intersection tree $\cT_{\cap,d}$ drawn from $GW^{(\lambda s)}_d$. Let $\cT_d$ be the shuffling of a $(\lambda,s,d)$ augmentation of $\cT_{\cap,d}$, and $\cT'_d$ be similarly obtained via a shuffling and an augmentation of $\cT_{\cap,d}$ independently performed. 
\end{definition}

\hbox{}

By the previous statement on the neighborhood of some vertex in the union graph $A_\cup$, one readily obtains the following
\begin{lemma}\label{lem:matched}
Consider two vertices $i$, $u$ of $A_1$, $A_2$ respectively, such that $\pi^*(i)=u$. Consider their respective graph neighborhoods $\cB_1(i,d)$, $\cB_2(u,d)$ up to distance $d$. When these neighborhoods are trees, let $i$ and $u$ be their respective roots, and use {\em shuffling} as just defined to label their non-root vertices, independently for $i$ and for $u$. Then for $d\le \kappa \ln n$, where 
$\kappa\ln(\lambda(2-s))<1/2$,  the corresponding neighborhoods are with high probability trees the pair of which, thus labeled, converges in total variation to $\cP^{(\lambda,s)}_d$.
\end{lemma}

\hbox{}

In contrast, if we instead pick two arbitrary vertices $i$, $u$ from each graph without enforcing the constraint $\pi^*(i)=u$ we have the following
\begin{lemma}\label{lem:unmatched}
Consider two vertices  $i$, $u$ each picked independently, uniformly at random from $A_1$and $A_2$ respectively. Consider their respective graph neighborhoods $\cB_1(i,d)$, $\cB_2(u,d)$ up to distance $d$. When these neighborhoods are trees, let $i$ and $u$ be their respective roots, and use {\em shuffling} to label their non-root vertices, independently for $i$ and $u$. Then for $d\le \kappa\ln n$, where $\kappa\ln(\lambda(2-s))<1/2$, the corresponding neighborhoods are with high probability trees, the pair of which converges in total variation to $\cP^{(\lambda)}_d:=GW^{(\lambda)}_d\otimes GW^{(\lambda)}_d$.
\end{lemma}

\hbox{}

These lemmas suggest as a design principle for graph alignment to match vertex $i$ in $A_1$ to vertex $u$ in $A_2$ only when the pair of neighborhoods $(\cB_1(i,d),\cB_2(u,d))$ has a far larger likelihood under $\cP^{(\lambda,s)}_d$ than under $\cP^{(\lambda)}_d$. 

\hbox{}

This idea almost works, but needs some refinements. To describe these, we need the following
\begin{definition}
For a vertex $i$ in graph $G$, and neighbor $j$ of $i$ in $G$, we define the oriented neighborhood $\cB_G(j\to i,d)$ of $j$ away from $i$ up to distance $d$  as the collection of vertices $k$ at distance at most $d$ from $j$ in the graph $G$ after removal of its edge $(i,j)$. We also denote by $\cB_G(j\to i,d)$ the subgraph of $G$ induced by this directed neighborhood.                                      
\end{definition}

\hbox{}

The algorithm is then defined as follows. 
\begin{definition}
Consider a pair of vertices $i$, $u$ of $A_1$, $A_2$ respectively. Assume that their respective neighborhoods $\cB_1(i,d+1)$, $\cB_2(u,d+1)$ are trees, where $d\le \kappa \ln n$ and $\kappa \ln(\lambda(2-s))<1/2$. 

Algorithm MP-Align matches $i$ to $u$ if $i$ (respectively $u$) admits three neighbors $j_1,j_2,j_3$ (respectively, $v_1,v_2,v_3$) such that for each pair $(j_k,v_k)$, $k=1,2,3$, one has
$$
\ln\left( L^{(\lambda,s)}_d(\cB_1(j_k\to i,d),\cB_2(v_k\to u,d))\right)\ge \tau (\lambda s)^d
$$
for some suitable constant threshold $\tau=\Theta(1)$, where for any two trees $(t,t')\in \cY_d$, we let $L^{(\lambda,s)}_d(t,t')$ denote the likelihood ratio
$$
L^{(\lambda,s)}_d(t,t'):=\frac{\cP^{(\lambda,s)}_d(t,t')}{\cP^{(\lambda)}_d(t,t')}\cdot
$$
\end{definition}

\hbox{}

Before we state the main theorem of this section, some remarks are in order. 
\begin{remark}
Although not obvious at first sight, MP-Align can be implemented as a message-passing algorithm. Indeed, assume that the roots of trees $t$, $t'\in \cY_d$ have respectively $c$ and $c'$ children. For $i\in[c]$, let $t_i$ denote the subtree of $t$ rooted at the $i$-th child of $t$'s root. Then, denoting by $\cS(k,c)$ the number of injections from $[k]$ to $[c]$, it follows from the definitions of $\cP^{\lambda,s}_d$, $\cP^{\lambda}_d$ that
$$
L^{(\lambda,s)}_d(t,t')=\sum_{k=0}^{c\wedge c'}\psi(k,c,c')\sum_{\sigma\in \cS(k,c),\sigma'\in \cS(k,c')}\prod_{u=1}^k L^{(\lambda,s)}_{d-1}(t_{\sigma(u)},t'_{\sigma'(u)}),
$$
where 
$$
\psi(k,c,c'):=\frac{\cp_{\lambda s}(k)\cp_{\lambda(1-s)}(c_k)\cp_{\lambda(1-s)}(c'-k)}{\cp_\lambda(c)\cp_\lambda(c')}\cdot
$$
This readily gives an iterative method for computing $L^{(\lambda,s)}_\delta(\cB_1(j\to i,d),\cB_2(v\to u,d))$ for any two oriented edges $(j\to i)$, $(v\to u)$ of $A_1$, $A_2$ respectively and for any $\delta\in [d]$. Initializing with $L^{(\lambda,s)}_0\equiv 1$, $L^{(\lambda,s)}_\delta(\cB_1(j\to i,d),\cB_2(v\to u,d))$ is obtained as a sum of a product of likelihood ratios $L^{(\lambda,s)}_{\delta-1}(\cB_1(k\to j,d),\cB_2(w\to v,d))$ for $k$ a neighbor of $j$ distinct from $i$ and $w$ a neighbor of $v$ distinct from $u$.

\hbox{}

The recursive operation can then be viewed as ``passing messages'' from pairs of oriented edges $(k\to j),(w\to v)$ to neighboring pairs of oriented edges $(j\to i),(v\to u)$. Indeed, this message-passing approach to computing likelihood ratios has been implemented and numerically tested in Piccioli et al. \cite{piccioli2022aligning} and in \cite{luca}.  In this sense it is reminiscent of {\bf belief propagation}, a celebrated algorithm and powerful analytical tool in the study of graphical models. See Jordan and Wainwright \cite{jordan} for an introduction to belief propagation, and Decelle et al. \cite{decelle} for its application as an analytical tool from statistical physics to high-dimensional statistics. More precisely, belief propagation is used in \cite{decelle} to predict that the threshold for computational hardness of performing {\em community detection}, or equivalently clustering, in a random graph known as the Stochastic Block Model, occurs precisely at the so-called Kesten-Stigum threshold from branching process theory.
%[citer Mezard-Montanari?]
\end{remark}

\hbox{}

\begin{remark}
MP-Align is a polynomial-time algorithm. Indeed, the number of operations performed in order to update some $L^{(\lambda,s)}_\delta(\cB_1(j\to i,d),\cB_2(v\to u,d))$ is upper-bounded by $c^2 c! c'!$, where $c$, $c'$ are the degrees of vertices $j$ and $v$ respectively. For sparse Erd\H{o}s-Rényi random graphs, it holds with high probability that the largest node degree is $(1+o(1))\ln n/\ln (\ln n)$. Thus 
$$
c!\le \exp((1+o(1))\ln(\ln n)\times \ln n/\ln (\ln n)))=\exp((1+o(1))\ln n).
$$
Thus for $d=O(\ln n)$, the number of operations needed to compute $L_d^{\lambda,s}$ is $O(n^{2+o(1)})$, hence the polynomial complexity of MP-Align.
\end{remark}

\hbox{}

We now need the following the following definition, which was suggested to us by Jiaming Xu \cite{jiaming_private}:
\begin{definition}      
For integer $d\ge 0$, we say that an algorithm for graph alignment is $d$-local if it returns a collection $\cC$ of pairs of vertices $(i,u)$ from graphs $A_1$, $A_2$ obtained by letting
$$
(i,u)\in \cC \hbox { if and only if }f(\cB_1(i,d),\cB_2(u,d))=1,
$$
where $f$ is a function from pairs of labeled graphs with a distinguished ``root'' vertex to $\{0,1\}$.
\end{definition}

\hbox{}

Note that this definition of a local algorithm differs from the notion of a {\em local search} algorithm mentioned in the Introduction.

\hbox{}

Clearly, MP-Align is a $d$-local algorithm in the sense just defined, where we enforce $d\le \kappa \ln n$ for $\kappa$ such that $\kappa \ln(\lambda(2-s))<1/2$. 
We are now ready to state the following
\begin{theorem}\label{thm:1}
Let $KL_d(\lambda,s)$ denote the Kullback-Leibler divergence from $\cP^{(\lambda,s)}_d$ to $\cP^{(\lambda)}_d$,
$$
KL_d(\lambda,s):=\sum_{(t,t')\in \cY_d^2}\cP^{(\lambda,s)}_d(t,t') \ln\left(L^{(\lambda,s)}_d(t,t')\right).
$$
Assume that
\begin{equation}\label{eq:kl_div}
\lim_{d\to\infty} KL_d(\lambda,s)=+\infty.
\end{equation}
Then for suitable choice of threshold parameter $\tau=\Theta(1)$, for $d=\Theta(\ln n)$, $d\le \kappa \ln n$ where $\kappa\ln(\lambda(2-s))<1/2$, with high probability MP-Align provides $\Omega(n)$ correct vertex matches, i.e. pairs $(i,u)$ such that $u=\pi^*(i)$, and returns only $o(n)$ incorrectly matched pairs. 

\hbox{}

If instead \eqref{eq:kl_div} does not hold, then any $d$-local algorithm which provides  $\Omega(n)$ correct matches with high probability must also provide $\Omega(n^2)$ incorrect matches with high probability.  
\end{theorem}

\hbox{}

For detailed proofs of this result and the other statements in this Section we refer the reader to \cite{luca}. A sketch of the argument is as follows.

\hbox{}

For the negative part, assume some function $f$ is such that with high probability,
$$
\sum_{i\in [n]} f(\cB_1(i,d),\cB_2(\pi^*(i),d))=\Omega(n).
$$
Then by Lemma \ref{lem:matched}, necessarily for $(t,t')\sim \cP^{(\lambda,s)}_d$, one must have
$$
x:=\cP^{(\lambda,s)}_d(f(t,t')=1)=\Omega(1).
$$
Let now $y:=\cP^{(\lambda)}_d(f(t,t')=1)$ denote the probability of the same event under distribution $\cP^{(\lambda)}_d$. Since we are under the assumption that \eqref{eq:kl_div} fails, the Kullback-Leibler divergence $KL_d(\lambda,s)$ is upper-bounded by some constant $C$. Classical monotonicity properties of the Kullback-Leibler divergence then imply that
$$
x\ln(x/y)+(1-x)\ln((1-x)/(1-y))\le C,
$$
and elementary inequalities allow to deduce
$$
y \ge x \exp(-C-1/e),
$$
hence $y=\Omega(1)$. 

\hbox{}

The number of incorrectly matched pairs reads
$$
\sum_{i,j\in[n],\pi^*(i)\ne j} f(\cB_1(i,d),\cB_2(j,d)).
$$
The expectation of this quantity is, by Lemma \ref{lem:unmatched}, $(1+o(1))n^2 y$. Arguing as in the proof of Lemma \ref{lem:unmatched}, one also gets that the second moment of the previous sum equals $(1+o(1))n^4 y^2$. Bienaymé-Tchebitchev inequality then ensures that with high probability, the number of incorrect matches is $(1+o(1))n^2 y=\Omega(n^2)$, thus obtaining the negative part of the Theorem.

\hbox{}

The positive part of Theorem \ref{thm:1} is shown by expressing the likelihood ratio $L_d^{(\lambda,s)}(t,t')$, for $(t,t')$ drawn according to $\cP^{(\lambda,s)}_d$ in a factorized form with factors $L_h^{(\lambda,s)}(t_i,t'_i)$ for $i$ in range $[\Theta((\lambda s)^d)]$, where the pairs $(t_i,t'_i)$ are drawn according to $\cP^{(\lambda,s)}_h$. Taking logarithms, the law of large numbers ensures that 
$$
\ln\left(\prod_{i=1}^{\Theta((\lambda s)^d)}L_h^{(\lambda,s)}(t_i,t'_i)\right)\approx \Theta((\lambda s)^d) KL_h(\lambda,s),
$$
and this is the leading contribution to $\ln(L^{(\lambda,s)}_d(t,t'))$ for sufficiently large $h$, the latter point following from the assumption \eqref{eq:kl_div} that $\lim_{h\to\infty}KL_h(\lambda,s)=+\infty$.
%%%
A direct consequence of Theorem \ref{thm:1} is the following
\begin{corollary}
MP-Align succeeds for $s>s^*(\lambda)$ where the threshold $s^*(\lambda)$ is defined via
\begin{equation}\label{eq:s_star}
s^*(\lambda):=\inf\{s\in [0,1]:\lim_{d\to\infty}KL_d(\lambda,s)=+\infty\}.
\end{equation}
\end{corollary}
 %%%
%%%%
%%%
\section{Tree correlation testing: Gaussian limits and Otter's constant}\label{sec:3}.
%%%%
%%%%
%%%
Let us introduce the following 
\begin{definition} Given two series of probability distributions $\{p_{0,d}\}_{d\ge 0}$, $\{p_{1,d}\}_{d\ge 0}$ on probability spaces that may depend on $d$, an associated series of tests $\{F_d\}$ taking values in $\{0,1\}$ is said to be one-sided if 
$$
\liminf_{d\to\infty} p_{1,d}(F_d=1)>0\hbox{ and }\lim_{d\to\infty}p_{0,d}(F_d=1)=0.
$$                                                         
\end{definition}
We then have the following characterization for existence of one-sided tests:
\begin{proposition}\label{prop:1}
Assume there is a one-sided family of tests for $\{p_{0,d}\}_{d\ge 0}$ against $\{p_{1,d}\}_{d\ge 0}$. Then necessarily, the Kullback-Leibler divergence $KL(p_{1,d}||p_{0,d})$ goes to infinity as $d\to\infty$. 

\hbox{}

Assume that there is a single probability space $(\Omega,\cF)$ endowed with two probability distributions $p_{1}$, $p_0$  and a filtration $\{\cF_d\}_{d\ge 0}$. Let $p_{1,d}$ (respectively, $p_{0,d}$) be the restriction of $p_1$ (respectively of $p_0$) to $\cF_d$. If the Kullback-Leibler divergence $KL(p_{1,d}||p_{0,d})$ goes to infinity as $d\to\infty$, then one-sided testability holds between $\{p_{1,d}\}_{d\ge 0}$ and $\{p_{0,d}\}_{d\ge 0}$. 
\end{proposition}  

\hbox{}
     
The first statement easily follows from the expression of the Kullback-Leibler divergence between two distributions $p,q$  as the supremum over finite measurable partitions of the Kullback-Leibler divergence between the discrete distributions induced by $p$ and $q$ on this finite partition, while the second follows from Doob's martingale convergence theorem.

\hbox{}

In the present context, let $\Omega=\cY_\infty\times \cY_\infty$, where $\cY_\infty$ is the space of rooted labeled trees, with no restriction on their depth, and take for $\cF_d$ the sigma-field generated by restriction to depth $d$ of the two trees.    Taking $p_{0,d}=\cP^{(\lambda)}_d$ and $p_{1,d}=\cP^{(\lambda,s)}_d$, Proposition \ref{prop:1} applies and yields that Condition \eqref{eq:kl_div} is equivalent to one-sided testability between correlated pairs of trees with distribution $\cP^{(\lambda,s)}_d$ and independent pairs of trees with distribution $\cP^{(\lambda)}_d$. 

\hbox{}

To gain further understanding into the threshold $s^*(\lambda)$ in \eqref{eq:s_star}, a fruitful approach consists in studying its limit as $\lambda\to\infty$. This will be done by establishing convergence to limiting Gaussian distributions of trees $t,t'$ sampled from $\cP^{(\lambda,s)}_d$ as $\lambda\to\infty$. 

\hbox{}

To that end one needs to consider rooted, unlabeled trees, i.e. equivalence classes of rooted labeled trees, rather than rooted labeled trees as in the last Section. Let then $\cX_d$ denote the collection of rooted unlabeled trees of depth at most $d$. An element of $\cX_d$ can be seen as a measure on $\cX_{d-1}$ with weights in $\dN$, which is adequate for establishing convergence to limiting Gaussian distributions in the spirit of the central limit theorem, as we are now dealing with elements from some vector space $\dR^{\cX_{d-1}}$. 

\hbox{}

For two elements $(\ct,\ct')$ of $\cX_d$, denote by $\cQ^{(\lambda,s)}_d(\ct,\ct')$, $\cQ^{(\lambda)}_d(\ct,\ct')$ the probability weights that distributions $\cP^{(\lambda,s)}_d$, $\cP^{(\lambda)}_d$ respectively put on all pairs $(t,t')\in\cY_d$ such that $t\in\ct$, $t'\in \ct'$, where we write $t\in \ct$ to say that $t$ belongs to the equivalence class $\ct$. We abuse notation and still write $L^{(\lambda,s)}_d$ for the likelihood ratio of these two distributions:
$$
L^{(\lambda,s)}_d(\ct,\ct')=\frac{\cQ^{(\lambda,s)}_d(\ct,\ct')}{\cQ^{(\lambda)}_d(\ct,\ct')},
$$
and also keep notation $GW^{(\lambda)}_d(\ct)$ to denote the probability that a depth $d$ Galton-Watson tree with offspring distribution $\cp_\lambda$ belongs to class $\ct$.

\hbox{}

The key step towards obtaining a tractable Gaussian limit is the following diagonalization of the likelihood ratio $L^{(\lambda,s)}_d$ on a suitable basis of orthogonal polynomials:
\begin{theorem}
For all $\lambda>0$, $d\in \dN$, $s\in[0,1]$, there exists a collection $\left\{\left(f^{(\lambda)}_{d,\beta}(x)\right)_{x\in \cX_d},\beta\in \cX_d
\right\}$ indexed by $\beta\in \cX_d$ of vectors $f^{(\lambda)}_{d,\beta}$ whose entries $f^{(\lambda)}_{d,\beta}(x)$ are polynomial in the coordinates of $x$, that do not depend on $s$, and such that
\begin{equation}\label{eq:diag}
\forall x,y\in \cX_d,\; L^{(\lambda,s)}_d(x,y)=\sum_{\beta\in \cX_d}s^{|\beta|-1}f^{(\lambda)}_{d,\beta}(x)f^{(\lambda)}_{d,\beta}(y),
\end{equation}
where for any tree $\beta\in \cX_d$, $|\beta|$ denotes its number of vertices. 
The polynomials $f^{(\lambda)}_{d,\beta}$ further satisfy the following properties:
\begin{eqnarray}
\forall x\in \cX_d, \; f^{(\lambda)}_{d,\bullet}(x)=1\quad \hbox{(value at trivial tree)},
\\
\forall \beta, \beta'\in \cX_d,\; \sum_{x\in\cX_d}GW^{(\lambda)}_d(x)f^{(\lambda)}_{d,\beta}(x)f^{(\lambda)}_{d,\beta'}(x)=\II_{\beta=\beta'}, \quad \hbox{(first orthogonality relation)},
\\
\forall x,x'\in \cX_d,\;\sum_{\beta\in\cX_d}f^{(\lambda)}_{d,\beta}(x)f^{(\lambda)}_{d,\beta}(x')=\frac{\II_{x=x'}}{GW^{(\lambda)}_d(x)},\quad\hbox{(second orthogonality relation)},
\end{eqnarray}
and finally for all $m\ge 2$, $d\ge 1$, $\beta^{(1)}=\{\beta^{(1)}_\gamma\}_{\gamma\in\cX_{d-1}}$,\ldots,$\beta^{(m)}=\{\beta^{(m)}_\gamma\}_{\gamma\in \cX_{d-1}}$, the mixed moments of order $m$ admit the following limits:
\begin{equation}\label{eq:mixed}
\lim_{\lambda\to\infty}\sum_{y\in\cX_d}GW^{(\lambda)}_d(y)f^{(\lambda)}_{d,\beta^{(1)}}(y)\cdots f^{(\lambda)}_{d,\beta^{(m)}}(y)=
\prod_{\gamma\in\cX_{d-1}}\sqrt{\prod_{i=1}^m\beta_\gamma^{(i)}!}\left[ x_1^{\beta^{(1)}_\gamma}\cdots x_m^{\beta^{(m)}_\gamma}\right]e^{\sum_{1\le i<j\le m}x_ix_j},
\end{equation}
where $[x^\alpha]G$ refers to the coefficient of monomial $x^\alpha$ in the power series expansion of some function $G$.
\end{theorem}

\hbox{}

The construction of polynomials $f_{d,\beta}^{(\lambda)}$ is done in \cite{ganassali2022statistical} by induction on $d$, as are shown the corresponding properties. 
\begin{remark} We now provide the intuition for the diagonalization formula \eqref{eq:diag} and how it separates dependency on $s$ via the terms $s^{|\beta|-1}$ from dependency on $\lambda$ via the terms $f^{(\lambda)}_\beta$. Define a Markov process $\{Z_t\}_{t\in \dR_+}$ on $\cX_d$ as follows. At unit rate, each edge of the current tree is suppressed as well as the sub-tree below that edge. At rate $\lambda$, for each $h=0,\ldots,d-1$, each vertex of the tree of depth $h$ obtains a new child and that child's genealogy tree up to depth $d-h-1$, the genealogy tree of that child being distributed according to $GW^{(\lambda)}_{d-h-1}$. 

\hbox{}

It is readily seen that this defines a reversible Markov process with invariant distribution $GW^{(\lambda)}_d$. 
Initializing the Markov process by taking $Z_0\sim GW^{(\lambda)}_d$, then a pair of trees $\ct,\ct'$ drawn from $\cQ^{(\lambda,s)}_d$ is obtained by letting $(\ct,\ct')=(Z_0,Z_{-\ln(s)})$. 

\hbox{}

Denote now by $H$ the infinitesimal generator of this Markov process. Formally (this is always true for Markov jump processes on finite state spaces, but requires a proof for countably infinite state spaces such as $\cX_d$), the matrix $P$ of transition probabilities $P_{xy}=\dP(Z_u=y|Z_0=x)$ in time $u$ is obtained as $P=\exp(u H)$. Let now $D$ denote the diagonal matrix $D_{x,y}=GW^{(\lambda)}_d(x)\II_{x=y}$. Reversibility entails that matrix $M:=D^{1/2}H D^{-1/2}$ is symmetric. One thus expects $M$ to be diagonalizable, with real spectrum $\{\nu_i\}_i\in \cI$ for some countable index set $\cI$, and associated collection of orthonormal eigenvectors $g_i$:
$$
M=D^{1/2}H D^{-1/2}=\sum_{i\in \cI}\nu_i g_i g_i^\top, 
$$
and thus the transition probabilities matrix $P$ reads
$$
P=D^{-1/2}\exp(u M) D^{1/2}=\sum_{i\in \cI} e^{u\nu_i} D^{-1/2}g_i g_i^\top D^{1/2}.
$$
Now, the likelihood ratio $L^{(\lambda,s)}_d(x,y)$ equals $(P D^{-1})_{x,y}$ for the specific choice $u=-\ln(s)$, so that
$$
L^{(\lambda,s)}_d=\sum_{i\in \cI} s^{-\nu_i}(D^{-1/2}g_i)(D^{-1/2}g_i)^\top.
$$
In this light, taking $\cX_d$ for indexing set $\cI$, we interpret  \eqref{eq:diag} as follows: the spectrum of infinitesimal generator $H$ is given by $\{-|\beta|\}_{\beta\in \cX_d}$, and the eigenvector $g_\beta$ of $M$ associated to $\nu_\beta=-|\beta|$ is given by $D^{1/2}f_\beta$.
\end{remark}

\hbox{}

%For any fixed $\beta_1,\ldots,\beta_k\in (\cX_d)^k$, the limits of higher order moments
%$$
%\lim_{\lambda\to\infty}\sum_{x\in \cX_d}GW^{\lambda}_d(x)f_{d,\beta_1}^\lambda(x)\times\cdots \times f_{d,\beta_k}^\lambda(x)
%$$
%are also computed by induction on $d$, and 
The expression \eqref{eq:mixed} for the limits of mixed moments implies the following asymptotic Gaussianity result:
\begin{theorem}
For $\cT\in \cX_{d+1}$, $\cT=\{N_{\tau}\}_{\tau\in \cX_d}$, consider its representation $\{\phi^{(\lambda)}_{d,\beta}(\cT)\}_{\beta\in \cX_d}$ in the dual basis $\{f_{d,\beta}^{(\lambda)}\}_{\beta\in \cX_d}$ defined by
$$
\forall \beta\in \cX_d,\;\phi^{(\lambda)}_{d,\beta}(\cT):=\sum_{\tau\in \cX_d}\frac{1}{\sqrt{\lambda}}f_{d,\beta}^{(\lambda)}(\tau)\left[N_\tau - GW^{(\lambda)}_d(\tau)\right].
$$
Then for a pair 
$(\cT=\{N_{\tau}\}_{\tau\in \cX_d}, \cT'=\{N'_\tau\}_{\tau\in \cX_d})$ drawn from $\cQ^{(\lambda,s)}_{d+1}$, the joint dual representations $\{\phi^{(\lambda)}_{d,\beta}(\cT),\phi^{(\lambda)}_{d,\beta}(\cT')\}_{\beta\in \cX_d}$ admit the following weak limit as $\lambda\to\infty$:
\begin{equation}
\lim_{\lambda\to\infty}\cQ^{(\lambda,s)}_{d+1}\left(\left\{\{\phi^{(\lambda)}_{d,\beta}(\cT),\phi^{(\lambda)}_{d,\beta}(\cT')\right\}_{\beta\in \cX_d}\in\cdot\right)=
%\stackrel[\beta\in \cX_d]{}{\otimes}
%\newcommand{\tens}[1]{%
%  \mathbin{\mathop{\otimes}\limits_{#1}}%
%}
\mathbin{\mathop{\otimes}\limits_{\beta\in \cX_d}}
\cN\left(0,\left(\begin{array}{cc} 1& s^{|\beta|}\\ s^{|\beta|}&1\end{array}\right)\right).
\end{equation}
\end{theorem}
Note that the distribution $\cQ^{(\lambda)}_d$ corresponds to the uncorrelated version ($s=0$) of $\cQ^{(\lambda,s)}_d$, and in the dual domain it thus also admits a Gaussian limit with the same structure as above, but with parameter $s=0$. Letting, for $k\in \dN$,
$$
J(k)=\left(\begin{array}{cc} 1& s^{k}\\ s^{k}&1\end{array}\right),
$$
and $I_2$ the $2\times 2$ identity matrix,
the Kullback-Leibler divergence between the limiting Gaussian distributions then evaluates to
\begin{equation}\label{eq:KL001}
\sum_{\beta\in \cX_d}KL\left(\cN\left(0,J(|\beta|)\right)||\cN\left(0,I_2\right)\right)=\sum_{\beta\in \cX_d}-\frac{1}{2}\ln(1-s^{2|\beta|})=\sum_{\beta\in \cX_d}\frac{1}{2}\sum_{j\ge 1} \frac{(s^j)^{2|\beta|}}{j}.
\end{equation}
Let $A_{d,n}$ denote the number of elements of $\cX_d$ with exactly $n$ vertices, and let for $x\in \dR_+$,
\begin{equation}
\phi_d(x):=\sum_{n\ge 1} A_{d,n} x^{n-1}.
\end{equation}
Using the correspondence between elements $\ct$ of $\cX_{d+1}$ and vectors $\{N_\tau\}_{\tau\in \cX_d}$ in $\dN^{\cX_d}$, it is shown in \cite{ganassali2022statistical} that
\begin{equation}\label{eq:KL002}
\phi_{d+1}(x)=\exp\left( \sum_{j\ge 1} \frac{x^j}{j} \phi_d(x^j)\right).
\end{equation}
It is further established  in \cite{ganassali2022statistical} that the limit as $\lambda\to\infty$ of $KL_d(\lambda,s)$ is the Kullback-Leibler divergence between the limiting distributions. We then obtain from  \eqref{eq:KL001}, \eqref{eq:KL002} the following
\begin{corollary}
Let $KL_d(\infty,s)$ denote the limit $\lim_{\lambda\to\infty}KL_d(\lambda,s)$. 
%Denote by $A_{n,d}$ the number of rooted unlabeled trees of depth at most $d$ and with $n$ vertices. 
One then has
\begin{equation}\label{eq:otter_1}
KL_d(\infty,s)=\frac{1}{2}\ln\left(\sum_{n\ge 1}A_{d,n}(s^2)^{n-1}\right).
\end{equation}
\end{corollary}

\hbox{}

Denote by $A_n$ the number of rooted unlabeled trees with $n$ vertices. The celebrated result of Otter \cite{otter1948number} establishes that for some constant $C>0$ one has
$$
A_n\sim \frac{C}{n^{3/2}}\alpha^{-n},
$$
where $\alpha=0.338...$ is now known as Otter's constant. This together with \eqref{eq:otter_1} in turn yields, since $A_{n,d}$ increases to $A_n$ as $d\to\infty$:
\begin{corollary}
If $s\le \sqrt{\alpha}$, then $\lim_{d\to\infty}KL_d(\infty,s)<+\infty$. 

If $s>\sqrt{\alpha}$, then $\lim_{d\to\infty}KL_d(\infty,s)=+\infty$.
\end{corollary}

\hbox{}

The relationship between Criterion \eqref{eq:kl_div} for success of Algorithm MP-Align and Otter's constant is then obtained by proving that the following interchange of limits holds:
$$
\lim_{d\to\infty} \lim_{\lambda\to\infty}KL_d(\lambda,s)=\infty\Leftrightarrow \lim_{\lambda\to\infty} \lim_{d\to\infty}KL_d(\lambda,s)=\infty.
$$
This finally yields the following result,  previously conjectured in Piccioli et al. \cite{piccioli2022aligning}:
\begin{theorem}
For $s\le \sqrt{\alpha}$, MP-Align fails. For $s>\sqrt{\alpha}$, there is $\lambda(s)$ such that, for $\lambda>\lambda(s)$, MP-Align succeeds. As a consequence, $\lim_{\lambda\to\infty} s^*(\lambda)=\sqrt{\alpha}$.
\end{theorem}

\begin{remark}        
In \cite{mao2023random}, Mao et al. propose an algorithm based on counting specific trees that they call chandeliers, and show that for correlated Erd\H{o}s-Rényi random graphs, with fixed $s>\sqrt{\alpha}$ and average degree $\lambda$ sufficiently large, their algorithm succeeds at partial alignment. The arguments differ from those in our present approach and provide a more direct justification for the appearance of Otter's constant. Also, their scheme is suited to non-sparse graphs ($\lambda\to\infty$) which MP-Align cannot handle.

\hbox{}

In contrast, our approach is specifically targeted at the sparse regime $\lambda=O(1)$. Also, Condition \eqref{eq:kl_div} for success of MP-Align coincides by Proposition \ref{prop:1} with one-sided testability between the correlated distributions $\cQ^{(\lambda,s)}_d$ and $\cQ^{(\lambda)}_d$, and by Theorem \ref{thm:1} with existence of some successful $d$-local algorithm for $d=O(\ln n)$. It is thus tempting to conjecture that Condition \eqref{eq:kl_div} precisely characterizes the onset of polynomial tractability in the sparse regime. This will be further discussed in Section \ref{sec:5}.
\end{remark}
\hbox{}

\section{Correlated Wigner matrices and a fast spectral method.}\label{sec:spectral}

Our aim in this Section is to identify a computational threshold for a spectral method called EIG1, to be defined shortly. EIG1 is the fastest among several such methods proposed in 
Feizi et al. \cite{feizi2019spectral}. We study it here in the specific context of correlated Gaussian random matrices 
$A_1$, $A_2$. Specifically, we are given two independent identically distributed matrices Gaussian Wigner $A_1$, $Z$. That is to say entries $\{Z(i,j)\}_{i\le j}$ are independent, with $Z(i,j)=Z(j,i)$ for all $i<j\in [n]$, $Z(i,i)$ i.i.d. with distribution $\cN(0,2/\sqrt{n})$ for $i\in [n]$, and $Z(i,j)$ i.i.d. with distribution $\cN(0,1/\sqrt{n})$ for $i<j\in [n]$.

\hbox{}

As before, $\pi^*$ is sampled uniformly at random from the set $\cS_n$ of permutations of $[n]$. We finally let 
$$
A'_2=A_1+\sigma Z,
$$ 
and 
$$
A_2(i,j)=A'_2(\pi^*(i),\pi^*(j)),\quad i\le j\in [n].
$$
The spectral method EIG1 proceeds as follows. Let $v_1$ (respectively, $v_2$) be a normalized eigenvector of $A_1$ (respectively, $A_2$) associated with its largest eigenvalue. Let $\pi_+$ (respectively, $\pi_-$) denote the permutation from $\cS_n$ that sorts the entries of eigenvectors $v_1$, $v_2$ in the same (respectively, in opposite) order, i.e. assuming that the entries of $v_1$, $v_2$ are sorted by permutations $\tau_1$, $\tau_2$, i.e.
$$
v_i(\tau_i(1))>v_i(\tau_i(2))>\cdots>v_i(\tau_i(n)),\quad i=1,2,
$$
then
$$
\forall j\in[n],\; \pi_+(\tau_1(j))=\tau_2(j)\hbox{ and }\pi_-(\tau_1(j))=\tau_2(n-j+1).
$$
Finally, the algorithm EIG1 returns the permutation $\hat\pi\in\{\pi_-,\pi_+\}$ chosen according to the criterion
$$
\hat\pi\in \hbox{argmax}_{\pi\in\{\pi_-,\pi_+\}}\left\{\sum_{i,j\in [n]}A_1(i,j)A_2(\pi(i),\pi(j))\right\}.
$$
EIG1's implementation requires computation of eigenvectors, which takes $\tilde{O}(n^2)$ operations, and sorting of the eigenvectors' entries, which takes $O(n\ln n)$ operations, thus an overall computational complexity of $\tilde{O}(n^2)$. As such it is faster to implement than competing proposals in the literature.

\hbox{}

Our main result for EIG1 is the following
\begin{theorem}\label{thm:spectr}
Consider the correlated Gaussian Wigner matrices $A_1$, $A_2$ with noise parameter $\sigma$. For any fixed $\epsilon>0$, if $\sigma\le n^{-7/6-\epsilon}$, then with high probability,
\begin{equation} \label{eq:eig_yes}
\hbox{ov}(\hat\pi,\pi^*)=1-o(1).
\end{equation}
If instead $\sigma\ge n^{-7/6+\epsilon}$, then with high probability,
\begin{equation}
\hbox{ov}(\hat\pi,\pi^*)=o(1).
\end{equation}
\end{theorem}

\hbox{}

The result is established in \cite{ganassali2022spectral}. The proof strategy consists in controlling the perturbation of the leading eigenvector of $A'_2$ that results from the addition of noise matrix $\sigma Z$ to $A_1$. In particular one establishes the following
\begin{proposition}\label{prop:22}
Assume that for some fixed $\alpha>0$, $\sigma\le n^{-1/2-\alpha}$. Then the leading eigenvector $v'_2$ of $A'_2$ can be written 
\begin{equation}
v'_2=(1+o_{\dP}(1))\left( v_1 +S\frac{1}{\|w\|}w\right),
\end{equation}
where $S$ is a random variable such that $S=(\sigma n^{1/6})^{1+o_\dP(1)}$, and $w$ is distributed as $\cN(0,I_n)$ and is independent of $v_1$.
\end{proposition}

\hbox{}

Proposition \ref{prop:22} essentially relies on writing the eigenvalue-eigenvector equation for eigenvector $v'_2$ of $A'_2$ in the basis of the eigenvectors of $A_1$. The evolution of the leading eigenvector $v'_2$ as $\sigma$ increases has been studied for instance by Allez and Bouchaud \cite{Allez_2014}, who establish stochastic differential equations for the evolution of the coordinates of a specific eigenvector of $A'_2$ in the eigenbasis of  $A_1$. 

\hbox{}

In order to obtain Proposition \ref{prop:22}, we need quantitative estimates on the magnitude of the perturbation. We obtain these by approximating the perturbed eigenvector via an iterative Picard scheme, whose convergence we analyze by exploiting fine rigidity estimates for the eigenvalues of Wigner matrices established in Erd\H{o}s, Yau and Yin \cite{ERDOS20121435}. 
%This analysis is carried through for $\sigma\le n^{-1/2-\epsilon}$ for some $\epsilon>0$, under which Condition 

\hbox{}

Equipped with Proposition \ref{prop:22}, the performance of EIG1 is reduced to the study of a toy problem, that consists in studying how the ranking of entries of a vector $x$ whose entries are i.i.d. with distribution $\cN(0,1/n)$ is affected by adding to each entry $x(i)$ a perturbation $S w(i)$ where $S=(\sigma n^{1/6})^{1+o_\dP(1)}$ and the $w(i)$ are i.i.d. $\cN(0,1)$, and independent of the $x(i)$. For the analysis of this toy problem, the proof of Proposition \ref{prop:22} and how they combine to yield Theorem \ref{thm:spectr}, we refer the reader to \cite{ganassali2022spectral}.
%propose several spectral methods for graph alignment, among which the simplest one

\hbox{}

\section{Correlated Wigner matrices and convex relaxation.}\label{sec:4}
In this Section we consider convex relaxations for the estimation of the unknown permutation $\pi^*$, placing ourselves again in the context of correlated Gaussian random matrices $A_1$, $A_2$. %Specifically, we are given two independent identically distributed matrices Gaussian Wigner $A_1$, $Z$. That is to say entries $\{Z(i,j)\}_{i\le j}$ are independent, with $Z(i,j)=Z(j,i)$ for all $i<j\in [n]$, $Z(i,i)$ i.i.d. with distribution $\cN(0,2/\sqrt{n})$ for $i\in [n]$, and $Z(i,j)$ i.i.d. with distribution $\cN(0,1/\sqrt{n})$ for $i<j\in [n]$.

\hbox{}

The algorithm considered in this Section returns matrix $X^*$, obtained as the minimizer of $X\to \|X A_1 -A_2 X\|_F$ over the convex set $\cB_n$, where $\cB_n$ is the set of doubly stochastic matrices, also known as the Birkhoff polytope:
%\sv{ Matrix $D$ notation seems to be inconsistent with the $M/X^\star$ notation below.}
\begin{equation}\label{eq:birkhoff}
X^*\in \hbox{argmin}_{X \in \cB_n}\|X A_1-A_2 X\|_F.
\end{equation}
%A popular choice for the convex set $\cC_n$ consists in letting $\cC_n=\cB_n$, that is the set of doubly stochastic matrices, also known as the Birkhoff polytope. By the Birkhoff-von Neumann theorem, $\cB_n$ is the convex hull of permutation matrices. 
%%
%%
Other convex relaxations are possible: for instance Araya and Tyagi \cite{araya2024graph} consider the above minimization over $M\in \Delta_n$, defined as the simplex in $\dR^{n^2}$:
$$
\Delta_n:=\left\{X\in \dR_+^{n^2}: \sum_{i,j\in [n]}X_{ij}=n\right\}.
$$
Surprisingly few theoretical results are available about the performance of such convex relaxations, despite their empirical success. Indeed, Araya and Tyagi \cite{araya2024graph} considered the noise-free case where $\sigma=0$, for which they show that for both relaxations to $\cB_n$ and to $\Delta_n$, the unique minimizer $X^*$ is given by permutation matrix $X^*=\Pi^*$ associated with $\pi^*$, for which the objective function $\|MA_1-A_2 M\|_F$ trivially equals zero. 

\hbox{}

The following result is established in Varma, Waldspurger and Massoulié \cite{sushil25}:
\begin{theorem}\label{thm:relax}
Consider $X^*$ that solves minimization problem \eqref{eq:birkhoff} over the Birkhoff polytope, for $A_1$, $A_2$ two correlated Gaussian Wigner matrices with noise parameter $\sigma\le 1$. 
Then for any fixed $\epsilon>0$, $X^*$ verifies with high probability:
\begin{eqnarray}
\sigma\le n^{-1-\epsilon}\Rightarrow \|X^*-\Pi^*\|^2_F=o(n) \label{eq:birkhoff_yes}\\
\sigma\ge n^{-1/2+\epsilon}\Rightarrow \|X^*-\Pi^*\|_F^2=\Omega(n). \label{eq:birkhoff_no}
\end{eqnarray}
\end{theorem}

\begin{remark}
In the first case \eqref{eq:birkhoff_yes}, property $\|X^*-\Pi^*\|^2_F=o(n)$ readily guarantees that a very crude post-processing of $X^*$, for instance letting
\begin{equation}\label{eq:post}
\forall i \in [n], \; \hat\pi(i)\in \hbox{argmax}_{j\in [n]}X^*_{ij}, 
\end{equation}
achieves overlap $\hbox{ov}(\hat \pi,\pi^*)$ of $1-o(1)$, i.e., it correctly recovers all but a vanishing fraction of the entries of $\pi^*$. Indeed, for any $i$ such that $\hat\pi(i)\ne \pi^*(i)$, it must be that $X^*_{i,\pi^*(i)}\le 1/2$. Thus 
$$
\|X^*-\Pi^*\|^2_F\ge \frac{1}{4}\sum_{i\in [n]}\II_{\hat\pi(i)\ne\pi^*(i)},
$$
hence the lower bound $1-o(1)$ on the overlap $\hbox{ov}(\hat\pi,\pi^*)$ when $\|X^*-\Pi^*\|^2_F=o(n)$. Since \eqref{eq:birkhoff_yes} is less stringent than Condition $\sigma\le n^{-7/6}$ needed for the success of EIG1 in \eqref{eq:eig_yes}, it follows that the convex relaxation has a better robustness to noise than EIG1, at the expense of higher computational complexity. 

\hbox{}

In contrast, under Condition \eqref{eq:birkhoff_no}, Property $\|X^*-\Pi^*\|_F^2=\Omega(n)$ is compatible with an overlap  $\hbox{ov}(\hat\pi,\pi^*)$ equal to $o(1)$. 
\end{remark}

\hbox{}

We briefly sketch how \eqref{eq:birkhoff_no} is established. Thanks to invariance of Frobenius norm $\|\cdot\|_F$ and of $\cB_n$ by left- or right-multiplication by a permutation matrix, we can assume without loss of generality that $\pi^*$ is the identity. Let $J$ be the all-ones matrix. Since $n^{-1}J\in \cB_n$, by optimality of $X^*$ one has
$$
\|A_1X^*-X^*A_2\|_F^2  \le \frac{1}{n^2}\|A_1 J-J A_2\|_F^2\le \frac{3}{n^2}\left[\|A_1 J\|_F^2+\|JA_1\|_F^2+\sigma^2\|JZ\|_F^2\right].
$$
The three Frobenius norms in the right-hand side are all of the form 
$$\sum_{i,j\in [n]}Y_{ij}^2,
$$
for $Y_{ij}\sim \cN(0,1+1/n)$. With high probability, all these Gaussian random variables are in absolute value less than $n^{\epsilon/2}$. This together with $\sigma\le 1$ ensures that with high probability,
\begin{equation}\label{eq:birkhoff_xxx}
\|A_1X^*-X^*A_2\|_F^2  \le 9 n^{\epsilon}.
\end{equation}
A companion lower bound on $\|A_1X^*-X^*A_2\|_F^2$ is obtained in terms of $\|X^*-I\|_F^2$ by writing
$$
\begin{array}{ll}
\|A_1X^*-X^*A_2\|_F^2 &=\|A_1(I-X^*)-(I-X^*)A_2+\sigma Z\|_F^2
\\
&=\|A_1(I-X^*)-(I-X^*)A_2\|_F^2+\sigma^2\|Z\|_F^2-2\sigma\langle A_1Z-ZA_1,X^*\rangle -2\sigma^2\langle Z^2,I-X^*\rangle.
\end{array}
$$
The first term $\|A_1(I-X^*)-(I-X^*)A_2\|_F^2$ is non-negative, the second term $\sigma^2\|Z\|_F^2$ is with high probability larger than $\sigma^2n/2$. We thus obtain, using Cauchy-Schwarz inequality for the fourth term:
$$
\|A_1X^*-X^*A_2\|_F^2 \ge \frac{\sigma^2 n}{2}-2\sigma\max_{i\ne j}|(A_1Z-ZA_1)_{ij}|\sum_{i\ne j}X^*_{ij} -2\sigma^2\|Z^2\|_F\|I-X^*\|_F.
$$
It can be shown that with high probability, $\max_{i\ne j}|(A_1Z-ZA_1)_{ij}|=O(n^{\epsilon/2-1/2})$. Noting further that $\sum_{i\ne j}X^*_{ij}\le n$, this yields
$$
\begin{array}{ll}
\|A_1X^*-X^*A_2\|_F^2 &\ge \frac{\sigma^2 n}{2}-O(\sigma n^{\epsilon/2+1/2})-2\sigma^2\|Z^2\|_F\|I-X^*\|_F
\\
&\ge \frac{\sigma^2 n}{2}-O(\sigma n^{\epsilon/2+1/2}) -5\sigma^2 \sqrt{n} \|I-X^*\|_F,
\end{array}
$$
where we used that $\|Z^2\|_F\le \sqrt{n}\|Z^2\|_2\le \sqrt{n}\|Z\|^2_2$ and that $\|Z\|_2=2+o(1)$ with high probability. By combining \eqref{eq:birkhoff_xxx} with this last display, and using that $\sigma \ge n^{-1/2+\epsilon}$, one gets
$$
5\sigma^2 \sqrt{n} \|I-X^*\|_F\ge \Omega(\sigma^2 n),
$$
hence the announced lower bound \eqref{eq:birkhoff_no}. We refer the reader to \cite{sushil25} for the detailed proof of Theorem \ref{thm:relax}. 

%The positive part is established by considering a dual convex optimization problem, and building a careful feasible dual solution based on the spectral decomposition of matrix $A_1$. 

%We are not aware of other results on such convex relaxations for correlated Gaussian Wigner matrices except for \cite{araya2024graph} and \cite{sushil25} just mentioned. We now state the main result of this Section, a tightening of the positive result in \eqref{eq:birkhoff}:
%\begin{theorem}\label{thm:relax}
%For correlated GOE matrices as above, if noise parameter $\sigma$ satisfies 
%\begin{equation}         
%\sigma=o\left( \frac{n^{-1/2}}{\ln(n)^4}\right),
%\end{equation}
%one then has with high probability, for both the relaxation to the Birkhoff polytope $\cB_n$ and to the constrained hypercube $\cH_n$ defined in \eqref{eq:hypercube}
%\begin{equation}\label{eq:borne_fr}
%\|X^*-\Pi^*\|^2_F=O(\sigma n^{3/2}\ln(n)^4)=o(n).
%\end{equation}
%\end{theorem}
%Before embarking on the proof, we remark that together with the negative result of \cite{sushil25} recalled in \eqref{eq:birkhoff}, Theorem \ref{thm:relax} locates precisely the threshold at which $\|X^*-\Pi^*\|_F^2$ goes from $o(n)$ to $\Omega(n)$ at $\sigma=n^{-1/2}$. 
%
\begin{remark} Post-processing \eqref{eq:post} is crude, and a more refined post-processing could be proposed. In particular one could let $\hat\pi$ be chosen as the permutation in $\cS_n$ that solves 
\begin{equation}\label{eq:postprocess}
\max_{\pi\in \cS_n}\sum_{i\in [n]}X^*_{i \pi(i)}.
\end{equation}
This is known as the Linear Assignment Problem (LAP), or maximum weight bipartite matching problem. Empirically LAP appears to achieve overlap $1-o(1)$ for $\sigma$ well above $n^{-1/2}$, potentially as large as $\Omega(1/\ln n)$, but there is so far no theoretical understanding of LAP's performance in the present setup.
\end{remark}

\hbox{}

\section{Perspectives.}\label{sec:5}

%No rigorous results of this kind are available to date.

Theorem \ref{thm:relax} from Section \ref{sec:4} predicts a threshold between $\sigma=n^{-1}$ and $\sigma=n^{-1/2}$ for the squared Frobenius norm $\|X^*-\Pi^*\|_F^2$ to go from $o(n)$ to $\Omega(n)$, a change which may compromise the accuracy of estimating $\pi^*$ by the simple post-processing of $X^*$ previously described by \eqref{eq:postprocess}. We are currently investigating the precise localization of this threshold, which we expect to be at $\sigma=n^{-1/2}$.

\hbox{}

As already mentioned, it is not known at which level of noise parameter $\sigma$ do more sophisticated methods such as LAP break down. Numerical experiments suggest that LAP produces a permutation $\hat\pi$ with overlap $1-o(1)$ for $\sigma$ potentially as large as $O(1/\ln n)$, and one promising direction consists in trying to get a precise characterization of the threshold for such methods.

\hbox{}

Section \ref{sec:spectral} focused on EIG1, one among several spectral methods proposed in Feizi et al. \cite{feizi2019spectral}. Another family of spectral methods described in \cite{feizi2019spectral} is as follows: construct an $n^2$ by $n^2$ matrix $M$, by letting
$$
M_{(i,k),(j,\ell)}=F(A_1(i,j),A_2(k,\ell))
$$
for some function $F:\dR^2\to \dR$. Proceed then to extract the eigenvector $\{x(i,k)\}_{i,k\in [n]}$ of $M$ corresponding to its largest eigenvalue. This can also be seen as a matrix, which may then be post-processed by using for instance LAP to obtain a permutation. The study of the spectrum of such random matrices is largely open. An exception is the case where $F(a,b)=a\times b$, which corresponds to $M=A_1\otimes A_2$ and for which the spectrum of $M$ is the tensor product of the spectra of $A_1$ and $A_2$, and the eigenvectors of $M$ are the tensor products of $A_1$ and $A_2$'s eigenvectors. Then essentially the method just described reduces to EIG1. However for more general functions $F$, the question of characterizing the spectral properties of the associated matrix $M$ for correlated GOE matrices $A_1$, $A_2$ constitutes an interesting challenge.

\hbox{}

%s mentioned in the Introduction, besides study of energy landscapes $\pi\to\dP(\pi^*=\pi|A_1,A_2)$, another way to show hardness of the estimation problem at hand consists in refuting the existence of estimators that are polynomials of the inputs $(A_1,A_2)$, with degree $D$ at most $O(\ln (n))$. Several results of this kind have been established under Condition $s<\sqrt{\alpha}$ for graph alignment. 

%First, Ding, Du and Li \cite{ding-du-li-23} have established low-degree hardness of the hypothesis testing version of graph alignment, that consists given a pair of graphs $(A_1,A_2)$ of deciding whether they are correlated Erd\H{o}s-Rényi graphs with parameters $(n,p,s)$ or alternatively whether they are independent $\cG(n,p)$ Erd\H{o}s-Rényi graphs when $s<\sqrt{\alpha}$.

%More recently, Li \cite{zhangsong-li-25} has shown, under an appropriate low degree conjecture, that positive overlap $\hbox{ov}(\hat \pi,\pi^*)=\Omega(1)$ cannot be achieved with low degree polynomials when $s<\sqrt{\alpha}$.

%Such hardness of partial graph alignment with low-degree polynomials when $s<\sqrt{\alpha}$, together with the fact that the threshold $s^*(\lambda)$ for success of MP-Align satisfies $\lim_{\lambda\to\infty}s^*(\lambda)=\sqrt{\alpha}$, tentatively suggests that in the sparse regime $\lambda=O(1)$, Condition \eqref{eq:kl_div} fails whenever partial graph alignment cannot be achieved with low-degree polynomials. We shall return to this point in Section \ref{sec:5}.

Sections \ref{sec:2} and \ref{sec:3} suggest another line of research, namely to try and obtain hardness results for graph alignment of sparse correlated Erd\H{o}s-Rényi graphs. Indeed Condition \eqref{eq:kl_div} is, as shown in Section \ref{sec:2}, the threshold for feasibility of partial alignment with $O(\ln n)$-local algorithms. It is tempting to guess that \eqref{eq:kl_div} is also the threshold for other notions of hardness, which would lend support to the belief that it is the threshold for polynomial feasibility of partial alignment. Progress along these lines would be of interest not just for graph alignment but also more broadly in high-dimensional statistics, by shedding new light on the relationships between the various notions of hardness.

\hbox{}

A first objective would aim to prove that partial alignment cannot be done with local search algorithms when \eqref{eq:kl_div} fails. More specifically, denote the posterior distribution of the sought permutation $\pi^*$ by $\dP_{post}(\tau):=\dP(\pi^*=\tau|A)$, where again $A=(A_1,A_2)$ is the observation. An impossibility result for local search algorithms would follow if under $\dP_{post}$, the overlap with $\pi^*$ satisfies a large deviations principle:
$$
\forall C\subset [0,1],\;\lim_{n\to \infty} \frac{1}{n}\ln(\dP_{post}(\{\tau\in\cS_n:\hbox{ov}(\tau,\pi^*)\in C))=-\inf_{x\in C}\phi(x),
$$ 
for some rate function $\phi$, and if furthermore 0 is a local minimum of $\phi$.  

\hbox{}

%could also be the threshold for the onset of polynomial tractability. 

%One first direction would consist in sh

%Showing that when it fails, either low-degree polynomials also fail, or the energy landscape $\pi\to \dP(\pi^*=\pi|A)$ becomes ill-behaved would be a worthy achievement, shedding new light on the relationships between the various notions of hardness considered in high-dimensional statistics. 

A second objective could consist in showing that partial alignment cannot be done with low-degree polynomial algorithms when \eqref{eq:kl_div} fails. 

\hbox{}

As mentioned in the Introduction, Ding, Du and Li \cite{ding-du-li-23} have established low-degree hardness of the hypothesis testing version of graph alignment, that consists given a pair of graphs $(A_1,A_2)$ of deciding whether they are correlated Erd\H{o}s-Rényi graphs with parameters $(n,p,s)$ or alternatively whether they are independent $\cG(n,p)$ Erd\H{o}s-Rényi graphs when $s<\sqrt{\alpha}$.

\hbox{}

More recently, Li \cite{zhangsong-li-25} has shown, under an appropriate low degree conjecture, that positive overlap $\hbox{ov}(\hat \pi,\pi^*)=\Omega(1)$ cannot be achieved with low degree polynomials when $s<\sqrt{\alpha}$. Li' approach \cite{zhangsong-li-25} reduces the task of proving LDP (low-degree polynomial) impossibility of estimation to proving LDP impossibility of an associated hypothesis testing task, modulo an adequate conjecture about LDP hardness. More precisely, the hypothesis testing problem considered in \cite{zhangsong-li-25}  is to distinguish between two independent Erd\H{o}s-Rényi graphs and two correlated Erd\H{o}s-Rényi graphs conditioned on the event $\pi^*(1)=1$. 

\hbox{}

The fact established in Section \ref{sec:2} that the threshold $s^*(\lambda)$ corresponding to  \eqref{eq:kl_div} converges to $\sqrt{\alpha}$ as $\lambda\to\infty$, together with the aforementioned results \cite{ding-du-li-23} and \cite{zhangsong-li-25} lends credit to the belief that, in the sparse regime $\lambda=O(1)$,  partial alignment is LPD-hard when \eqref{eq:kl_div}  fails. 

\hbox{}

Two possibilites can be envisioned to attack this problem. First, one could follow the approach of \cite{zhangsong-li-25} which reduces the problem to establishing LDP hardness of a hypothesis testing problem. Second, one could aim to establish directly LDP hardness of partial alignment, without relying on the conjecture used in \cite{zhangsong-li-25}. To that end, one could try the novel proof technique developed by Wein and Sohn \cite{10.1145/3717823.3718294} for showing LDP hardness of estimation, which they applied successfully to many problems (among which estimation of planted dense subgraphs and planted submatrices). The method of \cite{10.1145/3717823.3718294} in some sense leverages conditional independence properties of the models at hand, and may be applicable to graph alignment.

\hbox{}

Finally, the two models we considered, namely correlated Erd\H{o}s-Rényi graphs and correlated Gaussian Wigner matrices, both display an extreme amount of symmetries in the observations $A_1$, $A_2$. While they constitute a natural starting point for mathematical analysis, less symmetrical models may be closer to real-world data, and their analysis may reveal new phenomena of computational hardness. To date, few attempts at analyzing computational feasibility of graph alignment for less symmetrical models of graphs have been made. The main exceptions address correlated random graphs drawn from the Stochastic Block Model rather than Erd\H{o}s-Rényi graphs; see e.g. Li \cite{zhangsong-li-25} and Chen et al. \cite{CDGL}.
In this respect, Semerjian \cite{guilhem} provides an inspiring discussion of the role of symmetries in high dimensional statistics, focusing on how symmetries impact the feasibility of non-trivial estimation without any constraints on computational resources. Further investigations could consider the impact of symmetries on computational hardness of estimation, and aim to identify computational thresholds for less symmetrical yet tractable models. 

\hbox{}

%\appendix
\section*{Acknowledgements}. It is the author's pleasure to acknowledge stimulating interactions with his colleagues Luca Ganassali, Georgina Hall, Marc Lelarge, Jakob Maier, Guilhem Semerjian, Sushil Varma, Louis Vassaux, Irène Waldspurger and Jiaming Xu, who all contributed to the work discussed here. 

\hbox{}

\bibliographystyle{plain}
\bibliography{bib,references_jakob,references}

\begin{thebibliography}{10}

\bibitem{Allez_2014}
Romain Allez and Jean-Philippe Bouchaud.
\newblock Eigenvector dynamics under free addition.
\newblock {\em Random Matrices: Theory and Applications}, 03(03):1450010, July
  2014.

\bibitem{alon2016probabilistic}
Noga Alon and Joel~H Spencer.
\newblock {\em The probabilistic method}.
\newblock John Wiley \& Sons, 2016.

\bibitem{araya2024graph}
Ernesto Araya and Hemant Tyagi.
\newblock Graph matching via convex relaxation to the simplex.
\newblock {\em Foundations of Data Science}, pages 0--0, 2024.

\bibitem{franz-parisi}
Afonso~S. Bandeira, Ahmed~El Alaoui, Samuel~B. Hopkins, Tselil Schramm,
  Alexander~S. Wein, and Ilias Zadik.
\newblock The franz-parisi criterion and computational trade-offs in high
  dimensional statistics.
\newblock https://arxiv.org/abs/2205.09727, 2022.

\bibitem{bollobas2007phase}
B{\'e}la Bollob{\'a}s, Svante Janson, and Oliver Riordan.
\newblock The phase transition in inhomogeneous random graphs.
\newblock {\em Random Structures \& Algorithms}, 31(1):3--122, 2007.

\bibitem{CDGL}
Guanyi Chen, Jian Ding, Shuyang Gong, and Zhangsong Li.
\newblock A computational transition for detecting correlated stochastic block
  models by low-degree polynomials.
\newblock https://arxiv.org/abs/2409.00966, 2024.

\bibitem{conte2004thirty}
Donatello Conte, Pasquale Foggia, Carlo Sansone, and Mario Vento.
\newblock Thirty years of graph matching in pattern recognition.
\newblock {\em International journal of pattern recognition and artificial
  intelligence}, 18(03):265--298, 2004.

\bibitem{cullina-kyavash}
Daniel Cullina and Negar Kiyavash.
\newblock Exact alignment recovery for correlated erdős-rényi graphs.
\newblock https://arxiv.org/abs/1711.06783, 2017.

\bibitem{decelle}
Aurelien Decelle, Florent Krzakala, Cristopher Moore, and Lenka Zdeborová.
\newblock Inference and phase transitions in the detection of modules in sparse
  networks.
\newblock {\em Physical Review Letters}, 2011.

\bibitem{ding-hu}
Jian Ding and Hang Du.
\newblock Matching recovery threshold for correlated random graphs.
\newblock {\em The Annals of Statistics}, pages 1718--1743, 2023.

\bibitem{ding-du-li-23}
Jian Ding, Hang Du, and Zhangsong Li.
\newblock Low-degree hardness of detection for correlated erd\H{o}s-rényi
  graphs.
\newblock https://arxiv.org/abs/2212.13677, 2023.

\bibitem{DBLP:journals/focm/DingL25}
Jian Ding and Zhangsong Li.
\newblock A polynomial time iterative algorithm for matching gaussian matrices
  with non-vanishing correlation.
\newblock {\em Found. Comput. Math.}, 25(4):1287--1344, 2025.

\bibitem{ding2021efficient}
Jian Ding, Zongming Ma, Yihong Wu, and Jiaming Xu.
\newblock Efficient random graph matching via degree profiles.
\newblock {\em Probability Theory and Related Fields}, 179:29--115, 2021.

\bibitem{ERDOS20121435}
László Erdős, Horng-Tzer Yau, and Jun Yin.
\newblock Rigidity of eigenvalues of generalized wigner matrices.
\newblock {\em Advances in Mathematics}, 229(3):1435--1515, 2012.

\bibitem{DBLP:journals/focm/FanMWX23}
Zhou Fan, Cheng Mao, Yihong Wu, and Jiaming Xu.
\newblock Spectral graph matching and regularized quadratic relaxations {I}
  algorithm and gaussian analysis.
\newblock {\em Found. Comput. Math.}, 23(5):1511--1565, 2023.

\bibitem{fan2023spectral2}
Zhou Fan, Cheng Mao, Yihong Wu, and Jiaming Xu.
\newblock Spectral graph matching and regularized quadratic relaxations {II}:
  Erd{\H{o}}s-r{\'e}nyi graphs and universality.
\newblock {\em Foundations of Computational Mathematics}, 23(5):1567--1617,
  2023.

\bibitem{feizi2019spectral}
Soheil Feizi, Gerald Quon, Mariana Recamonde-Mendoza, Muriel Medard, Manolis
  Kellis, and Ali Jadbabaie.
\newblock Spectral alignment of graphs.
\newblock {\em IEEE Transactions on Network Science and Engineering},
  7(3):1182--1197, 2019.

\bibitem{DBLP:conf/msml/Ganassali21}
Luca Ganassali.
\newblock Sharp threshold for alignment of graph databases with gaussian
  weights.
\newblock In Joan Bruna, Jan~S. Hesthaven, and Lenka Zdeborov{\'{a}}, editors,
  {\em Mathematical and Scientific Machine Learning, 16-19 August 2021, Virtual
  Conference / Lausanne, Switzerland}, volume 145 of {\em Proceedings of
  Machine Learning Research}, pages 314--335. {PMLR}, 2021.

\bibitem{ganassali2022spectral}
Luca Ganassali, Marc Lelarge, and Laurent Massouli{\'e}.
\newblock Spectral alignment of correlated gaussian matrices.
\newblock {\em Advances in Applied Probability}, 54(1):279--310, 2022.

\bibitem{luca}
Luca Ganassali, Marc Lelarge, and Laurent Massoulié.
\newblock Correlation detection in trees for planted graph alignment.
\newblock {\em The Annals of Applied Probability}, 34(3), June 2024.

\bibitem{DBLP:conf/colt/GanassaliML21}
Luca Ganassali, Laurent Massouli{\'{e}}, and Marc Lelarge.
\newblock Impossibility of partial recovery in the graph alignment problem.
\newblock In Mikhail Belkin and Samory Kpotufe, editors, {\em Conference on
  Learning Theory, {COLT} 2021, 15-19 August 2021, Boulder, Colorado, {USA}},
  volume 134 of {\em Proceedings of Machine Learning Research}, pages
  2080--2102. {PMLR}, 2021.

\bibitem{ganassali2022statistical}
Luca Ganassali, Laurent Massoulié, and Guilhem Semerjian.
\newblock Statistical limits of correlation detection in trees.
\newblock {\em The Annals of Applied Probability}, 34(4), August 2024.

\bibitem{DBLP:journals/ior/HallM23}
Georgina Hall and Laurent Massouli{\'{e}}.
\newblock Partial recovery in the graph alignment problem.
\newblock {\em Oper. Res.}, 71(1):259--272, 2023.

\bibitem{thesis-hopkins}
Samuel Hopkins.
\newblock Statistical inference and the sum of squares method.
\newblock {\em PhD Thesis, Cornell University}, pages 1794--1827, 2018.

\bibitem{jordan}
Michael Jordan and Martin Wainwright.
\newblock {\em Graphical Models, Exponential Families, and Variational
  Inference}.
\newblock Now Publishers Inc, 2008.

\bibitem{zhangsong-li-25}
Zhangsong Li.
\newblock Algorithmic contiguity from low-degree conjecture and applications in
  correlated random graphs.
\newblock https://arxiv.org/abs/2502.09832, 2025.

\bibitem{jakob}
Jakob Maier and Laurent Massouli{\'e}.
\newblock Asymmetric graph alignment and the phase transition for asymmetric
  tree correlation testing.
\newblock https://arxiv.org/abs/2504.02299, 2025.

\bibitem{mao2023exact}
Cheng Mao, Mark Rudelson, and Konstantin Tikhomirov.
\newblock Exact matching of random graphs with constant correlation.
\newblock {\em Probability Theory and Related Fields}, 186(1):327--389, 2023.

\bibitem{mao2023random}
Cheng Mao, Yihong Wu, Jiaming Xu, and Sophie~H Yu.
\newblock Random graph matching at otter’s threshold via counting
  chandeliers.
\newblock In {\em Proceedings of the 55th Annual ACM Symposium on Theory of
  Computing}, pages 1345--1356, 2023.

\bibitem{narayanan2008robust}
Arvind Narayanan and Vitaly Shmatikov.
\newblock Robust de-anonymization of large sparse datasets.
\newblock In {\em 2008 IEEE Symposium on Security and Privacy (sp 2008)}, pages
  111--125. IEEE, 2008.

\bibitem{otter1948number}
Richard Otter.
\newblock The number of trees.
\newblock {\em Annals of Mathematics}, pages 583--599, 1948.

\bibitem{piccioli2022aligning}
Giovanni Piccioli, Guilhem Semerjian, Gabriele Sicuro, and Lenka Zdeborov{\'a}.
\newblock Aligning random graphs with a sub-tree similarity message-passing
  algorithm.
\newblock {\em Journal of Statistical Mechanics: Theory and Experiment},
  2022(6):063401, 2022.

\bibitem{guilhem}
Guilhem Semerjian.
\newblock Some observations on the ambivalent role of symmetries in bayesian
  inference problems.
\newblock {\em Comptes Rendus de l'Académie des Sciences}, 2025.

\bibitem{singh2008global}
Rohit Singh, Jinbo Xu, and Bonnie Berger.
\newblock Global alignment of multiple protein interaction networks with
  application to functional orthology detection.
\newblock {\em Proceedings of the National Academy of Sciences},
  105(35):12763--12768, 2008.

\bibitem{10.1145/3717823.3718294}
Youngtak Sohn and Alexander~S. Wein.
\newblock Sharp phase transitions in estimation with low-degree polynomials.
\newblock In {\em Proceedings of the 57th Annual ACM Symposium on Theory of
  Computing}, STOC '25, page 891–902, New York, NY, USA, 2025. Association
  for Computing Machinery.

\bibitem{sushil25}
Sushil Varma, Irène Waldspurger, and Laurent Massoulié.
\newblock Graph alignment via birkhoff relaxation.
\newblock https://arxiv.org/pdf/2503.05323, 2025.

\bibitem{Wein-arxiv}
Alexander~S. Wein.
\newblock Computational complexity of statistics: New insights from low-degree
  polynomials.
\newblock https://arxiv.org/abs/2506.10748, 2025.

\bibitem{jiaming_private}
Jiamig Xu.
\newblock Private communication, 2024.

\end{thebibliography}

%\appendix
%\end{thebibliography}

\end{document}